\documentclass[12pt,reqno]{amsart}
\usepackage[T1]{fontenc}
\usepackage{amsmath,amssymb,amsxtra,amscd}
\usepackage{amsthm,Transfer_plfgNotation}
\usepackage[all]{xy}

\newtheorem{theorem}{Theorem}[section]
\newtheorem{proposition}[theorem]{Proposition}
\newtheorem{lemma}[theorem]{Lemma}
\newtheorem{definition}[theorem]{Definition}

\newtheorem{corollary}[theorem]{Corollary}

\newtheorem*{question*}{Question}

\theoremstyle{definition}
\newtheorem{remark}[theorem]{Remark}

\begin{document}

\title{Retractive transfers and $p$-local finite groups}
\author{K\'ari Ragnarsson}
\date{October 19th, 2005}

 \address{Department of Mathematical Sciences, University of Aberdeen, Aberdeen \mbox{AB24 3UE}, United Kingdom}
\email{kari@maths.abdn.ac.uk}

\begin{abstract}
In this paper we explore the possibility of defining $p$-local
finite groups (\cite{BLO2}) in terms of transfer properties of their
classifying spaces. More precisely, we consider the question posed
by Haynes Miller, whether an equivalent theory can be obtained by
studying triples $(f,t,X)$, where $X$ is a $p$-complete, nilpotent
space with finite fundamental group, \mbox{$f\negmedspace: BS \to
X$} is a map from the classifying space of a finite $p$-group, and
$t$ is a stable retraction of $f$ satisfying Frobenius reciprocity
at the level of stable homotopy. We refer to $t$ as a
\emph{retractive transfer} of $f$ and to $(f,t,X)$ as a
\emph{retractive transfer triple over $S$}.

In the case where $S$ is elementary abelian, we answer this question
in the affirmative by showing that a retractive transfer triple
$(f,t,X)$ over $S$ does indeed induce a $p$-local finite group over
$S$ with $X$ as its classifying space.

Using results from \cite{KR:ClSpec} we show that the converse is
true for general finite $p$-groups. That is, for a $p$-local finite
group $\plfg$, the natural inclusion \mbox{$\theta \negmedspace: BS
\to X$} has a retractive transfer $t$, making $(\theta,t,\ClSp)$ a
retractive transfer triple over $S$. This also requires a proof,
obtained jointly with Ran Levi, that $\ClSp$ is a nilpotent space,
which is of independent interest.
\end{abstract}

\maketitle

\section*{Introduction} Defined by Broto-Levi-Oliver in
\cite{BLO2}, $p$-local finite groups are the culmination of a
program initiated by Puig \cite{Puig,Puig2} to find a formal
framework for the $p$-local structure of a finite group. To a finite
group $G$, one associates a fusion system (at a prime $p$)
consisting of all $p$-subgroups of $G$ and the homomorphims between
them induced by conjugation in $G$. Puig formalized fusion systems
and identified an important subclass of fusion systems, which we now
call saturated fusion systems. Fusion systems of finite groups are
contained in this class, but saturated fusion systems also arise in
other important contexts, most notably in modular representation
theory through Brauer subpairs of blocks of group algebras, and more
recently as Chevalley groups of $p$-compact groups \cite{BM}.

The fusion system of a group $G$ can be considered as an algebraic
interpretation of the $p$-local structure of the group. One can also
take a topological approach, and think of the $p$-local structure of
$G$ as being the $p$-completed classifying space $\pComp{BG}$. By
Bob Oliver's solution \cite{Ol1,Ol2} of the Martino-Priddy
conjecture \cite{MP}, these approaches are the same. That is, two
groups induce the same fusion system if and only if their
$p$-completed classifying spaces are homotopy equivalent. In fact
the fusion system can be recovered from the classifying space via a
homotopy theoretic construction, which is presented in Section
\ref{sub:Spaces-to-plfg}. Therefore, one can in some sense regard
\pComp{BG} as a classifying space of the fusion system. This
suggests that, more generally, each saturated fusion system may have
a unique classifying space.

A $p$-local finite group consists of a saturated fusion system and
an associated centric linking system, a category which contains just
enough information to construct a classifying space associated to
the fusion system. Thus one can think of a $p$-local finite group as
a saturated fusion system with a chosen classifying space.

The definition of $p$-local finite groups is rather complicated and
has the drawback that there is no straightforward concept of
morphisms between $p$-local finite groups, so they have not yet been
made to form a category in any sensible way. In this paper, we adopt
the approach used by Dwyer and Wilkerson for $p$-compact groups
\cite{DW}, and try to develop the theory of $p$-local finite groups
in terms of classifying spaces.

Specifically, we consider a homotopy monomorphism \mbox{$f
\negmedspace: BS \to X$} from the classifying space of a finite
$p$-group $S$ to a $p$-complete, nilpotent space $X$ with finite
fundamental group, that is endowed with a stable retraction $t$
satisfying Frobenius reciprocity at the level of stable homotopy. We
refer to such a triple $(f,t,X)$ as a \textit{retractive transfer
triple over $S$}. For a retractive transfer triple $(f,t,X)$, we ask
whether $X$ is the classifying space of a $p$-local finite group.
This question is addressed in Sections \ref{sec:Cohomology} and
\ref{sec:Homotopy}, where we answer the question in the affirmative
in the case where $S$ is elementary abelian (Theorem
\ref{thm:ElAb}). Conversely we ask whether a $p$-local finite group
$\plfg$ gives rise to a Frobenius transfer triple. This is indeed
the case for any finite $p$-group $S$, as we show in Section
\ref{sec:p-lfg-to-RTT}. This involves joint work with Ran Levi where
we show that the classifying space of a $p$-local finite group is a
torsion space and (since its fundamental group is a finite
$p$-group) consequently a nilpotent space.

The homotopy theory of classifying spaces of elementary abelian
$p$-groups was intensely studied in the nineties and is now well
understood through contributions by various authors, the most
important being Miller's solution of the Sullivan conjecture
\cite{Mil:Sullivan} and Lannes's $T$-functor technology \cite{Lan}.
Other contributions, that are related to the results and methods in
this paper are the work of Goerss--Smith-Zarati \cite{GSZ},
Harris--Kuhn \cite{HK}, Henn \cite{Henn}, Dwyer--Wilkerson,
\cite{DW2} and Henn--Lannes--Schwartz \cite{HLS}. The work in
Section \ref{sec:Homotopy} of this paper mimics the methods used by
Dwyer--Miller--Wilkerson in \cite{DMW}, replacing inclusions of
maximal tori in compact Lie groups with the inclusion of elementary
abelian $p$-groups in finite groups.

\textbf{Acknowledgements:} The majority of the work in this paper
was done as part of my Ph.D. thesis and I would like to thank my
thesis advisor Haynes Miller for suggesting this problem and for his
many helpful suggestions during the course of this project. I also
thank Ran Levi for his help on the material in Section
\ref{sub:Technical}. Finally, I thank the referee for suggesting the
more elegant approach of using the methods from \cite{DMW} in
Section \ref{sec:Homotopy} to construct an explicit map between
classifying spaces rather than using Wojtkowiak's obstruction theory
\cite{Wo1,Wo2} to establish its existence as I had originally
intended.

\section{A quick review of $p$-local finite groups}
In this section we give a brief overview of the theory of
$p$-local finite groups. Most of this material is found in
\cite{BLO2}. In this section and throughout the paper, $p$ is a
fixed prime.

\subsection{Some definitions and terminology}\label{sec:Def&Term}
We begin by recalling some terminology regarding $p$-local finite
groups.

\begin{definition} \label{def:FS}
A \emph{fusion system} $\F$ over a finite $p$-group $S$ is a
category, whose objects are the subgroups of $S$, and whose morphism
sets $\HomF{P}{Q}$ satisfy the following conditions:
\begin{itemize}
  \item[(a)] $\Hom{P}{Q}{S} \subseteq \HomF{P}{Q} \subseteq \Inj{P}{Q}{}$
for all $P,Q \leq S$.
  \item[(b)] Every morphism in $\F$ factors as
an isomorphism in $\F$ followed by an inclusion.
\end{itemize}
\end{definition}

Here $\Hom{P}{Q}{S}$ is the set of group homomorphims induced by
conjugation by elements in $S$.

Before stating the next definition, we need to introduce some
additional terminology and notation. We say that two subgroups $P,P'
\leq S$ are \emph{ $\F$-conjugate} if they are isomorphic in $\F$. A
subgroup $P \leq S$ is \emph{fully centralized in \F} if $\vert
C_S(P)\vert \geq \vert C_S(P')\vert$ for every $P' \leq S$ which is
$\F$-conjugate to $P$. Similarly $P$ is \emph{fully normalized in
\F} if $\vert N_S(P)\vert \geq \vert N_S(P')\vert$ for every $P'
\leq S$ which is $\F$-conjugate to $P$.

\begin{definition} \label{def:SatFS}
 A fusion system $\F$ over a $p$-group $S$ is \emph{saturated} if the
 following two conditions hold:
\begin{itemize}
  \item[(I)] If $P \leq S$ is fully normalized in \F, then $P$ is
also fully centralized, and $\Aut{P}{S}$ is a Sylow subgroup of
$\Aut{P}{\F}$.
  \item[(II)] If $P \leq S$ and \mbox{$\varphi
\in \HomF{P}{S}$} are such that $\varphi (P)$ is fully centralized,
then $\varphi$ extends to \mbox{$\bar \varphi \in
\HomF{N_{\varphi}}{S}$}, where
$$ N_{\varphi} = \N{\varphi}{P}{S}{g}.$$
\end{itemize}
\end{definition}

There is a class of subgroups of $S$ of special interest to us,
defined as follows.

\begin{definition}\label{def:Centric}
Let $\F$ be a fusion system over a $p$-group $S$. A subgroup
\mbox{$P \leq S$} is \emph{ $\F$-centric} if \mbox{$C_S(P') \leq
P'$} for every $P'$ that is $\F$-conjugate to $P$. Let $\F^{c}$
denote the full subcategory of $\F$ whose objects are the
$\F$-centric subgroups of $S$.
\end{definition}

\begin{remark}
The condition \mbox{$C_S(P') \leq P'$} in the previous definition
is equivalent to the condition \mbox{$C_S(P') = Z(P')$}.
\end{remark}

\begin{definition} \label{def:CLS}
Let $\F$ be a fusion system over the $p$-group $S$. A \emph{centric
linking system associated to} $\F$ is a category $\Link$ whose
objects are the $\F$-centric subgroups of $S$, together with a
functor
$$ \pi: \Link \rightarrow \F^{c},$$
and distinguished monomorphisms \mbox{$P \xrightarrow{\delta_P}
\Aut{P}{\Link},$} for each $\F$-centric subgroup $P \leq S$, which
satisfy the following conditions.

\begin{itemize}
  \item[(A)] The functor $\pi$ is the identity on objects and
surjective on morphisms. More precisely, for each pair of objects
$P,Q \in \Link,$ the centre $Z(P)$ acts freely on $\MorL{P}{Q}$ by
composition (upon identifying $Z(P)$ with \mbox{$\delta_P(Z(P))
\leq \Aut{P}{\Link}$}), and $\pi$ induces a bijection

\[
\begin{CD}
{\MorL{P}{Q}/Z(P)} @ > {\cong} >> {\HomF{P}{Q}.} \\
\end{CD}
\]

  \item[(B)] For each $\F$-centric subgroup $P \leq S$ and each $g
\in P$, $\pi$ sends $\delta_{P}(g) \in \Aut{P}{\Link}$ to $c_{g}
\in \Aut{P}{\F}$.

  \item[(C)] For each $f \in \MorL{P}{Q}$ and each $g \in P$, the
following square commutes in $\Link$:

\[
\begin{CD}
{P} @ > {f} >> {Q} \\ @ VV {\delta_{P}(g)} V @ VV
{\delta_{Q}(\pi(f)(g))} V \\ {P} @> {f} >> {Q.} \\
\end{CD}
\]
\end{itemize}
\end{definition}

We can now finally define our objects of study.
\begin{definition} \label{plfg}
A \emph{$p$-local finite group} is a triple $\plfg$, where $\F$ is a
saturated fusion system over a finite $p$-group $S$ and $\Link$ is a
centric linking system associated to $\F$. The \emph{classifying
space} of the $p$-local finite group is the $p$-completed geometric
realisation $\ClSp$.
\end{definition}

We remark that a $p$-local finite group comes equipped with a
natural inclusion
$$\theta\negmedspace: BS \longrightarrow \ClSp.$$

One of the main questions in the theory of $p$-local finite groups
concerns the existence and uniqueness of a centric linking system
associated to a given saturated fusion system. In Section 3 of
\cite{BLO2}, Broto-Levi-Oliver have developed an obstruction theory
to address this question.

\subsection{The fusion system of a group}\label{sec:Motivation} In
this section we will discuss the fusion system arising from a
Sylow subgroup inclusion \mbox{$S \leq G$}. This section serves as
motivation for the discussion in the previous section as well as
being of independent interest.

\begin{definition}\label{def:GroupFS}
Let $G$ be a finite group. The \emph{fusion system of $G$} is the
category $\F(G)$, whose objects are the $p$-subgroups of $G$ and
whose morphism sets are given by
$$\Hom{P}{Q}{\F(G)} = \Hom{P}{Q}{G}$$
for all $p$-subgroups \mbox{$P,Q \leq G$}.

For a $p$-subgroup $S \leq G$, the \emph{fusion system of $G$ over
$S$} is the full subcategory \mbox{$\F_S(G) \subseteq \F(G)$}, whose
objects are the subgroups of $S$.
\end{definition}

If $S$ is a Sylow subgroup of $G$, then the inclusion of $\F_S(G)$
in $\F(G)$ is an equivalence of categories, since every
$p$-subgroup of $G$ is conjugate to a subgroup of $S$.

\begin{proposition}\cite[Proposition 1.3]{BLO2}\label{prop:GroupFSSaturated}
Let $G$ be a finite group and let $S$ be a $p$-subgroup. Then the
fusion system $\F_S(G)$ of $G$ over $S$ is saturated if and only
if $S$ is a Sylow subgroup.
\end{proposition}

The centric linking system of a finite group was initially
introduced in \cite{BLO1} as a powerful tool to study homotopy
equivalences between $p$-completed classifying spaces of finite
groups. The $p$-centric subgroups of a finite group $G$ are the
$p$-subgroups \mbox{$P \leq G$} whose centre $Z(P)$ is a $p$-Sylow
subgroup of the centraliser $C_G(P)$. This notion of centricity is
equivalent to the one introduced in \ref{def:Centric} in the sense
that if \mbox{$S \leq G$} is a Sylow subgroup, then a subgroup $P
\leq S$ is $p$-centric if and only if it is $\F_S(G)$-centric.

For the following definition, we recall that if a group \mbox{$P
\leq G$} is $p$-centric, then one can write
$$C_G(P) = Z(P) \times C'_G(P),$$
where \mbox{$C'_G(P) \leq G$} has order prime to $p$. The notation
$C'_G(P)$ will be used in the definition. In addition, for
subgroups \mbox{$P,Q \leq G,$} we will let $N_G(P,Q)$ denote the
transporter
$$N_G(P,Q) = \{g \in G \vert gPg^{-1} \leq Q \}.$$

\begin{definition} \label{def:GroupCLS}
Let $G$ be a finite group. The \emph{centric linking system of $G$}
is the category $\Link(G)$, whose objects are the $p$-centric
subgroups of $G$ and whose morphism sets are given by
$$\Mor{P}{Q}{\Link(G)} = N_G(P,Q)/C'_G(P)$$
for all $p$-subgroups \mbox{$P,Q \leq G$}.

For a $p$-subgroup $S \leq G$, the \emph{centric linking system of
$G$ over $S$} is the full subcategory \mbox{$\Link_S(G) \subseteq
\Link(G)$} whose objects are the subgroups of $S$, that are
$p$-centric in $G$.
\end{definition}

In the case of a Sylow inclusion \mbox{$S \leq G$}, the centric
linking system $\Link_S(G)$ is a centric linking system associated
to the saturated fusion system $\F_S(G)$, and we have the following
proposition, which serves as a motivating example for the definition
of a $p$-local finite group.

\begin{proposition} \label{prop:Group-pLFG}
Let $S$ be a Sylow subgroup of a finite group $G$. Then the triple
$(S,\F_S(G),\Link_S(G))$ is a $p$-local finite group over $S$.
Furthermore, the natural map
$$\theta \negmedspace :BS \to \pComp{|\Link_S(G)|}$$
is equivalent to the $p$-completed inclusion
$$BS \to \pComp{BG}$$
as a space under $BS$.
\end{proposition}

\subsection{Homotopy theoretic constructions of fusion
systems}\label{sub:Spaces-to-plfg} \label{sec:HomotopyFS} In this
section we recall how a map \mbox{$f \negmedspace: BS \to X$}, from
the classifying space of a finite $p$-group $S$ to a space $X$,
induces a fusion system $\F_{S,f}(X)$ over $S$. In general this
fusion system is not saturated.

The following definition is motivated by the fact that two group
homomorphisms \mbox{$\varphi, \psi \negmedspace : G \to H$}
between finite groups are $H$-conjugate if and only if the induced
maps of classifying spaces are freely homotopic.

\begin{definition}\label{def:HomotopyFS}
For any space $X$, any $p$-group $S$, and any map \mbox{$f
\negmedspace: BS \to X$}, define $\F_{S,f}(X)$ to be the category
whose objects are the subgroups of $S$, and whose morphisms are
given by
$$\Hom{P}{Q}{\F_{S,f}} = \{\varphi \in \Inj{P}{Q}{} \mid f\vert_{BP} \simeq f\vert_{BQ}\circ B\varphi
\}$$ for each \mbox{$P,Q \leq S.$}
\end{definition}

It is easy to see that $\F_{S,f}$ is indeed a fusion system,
although it need not be saturated. In the case where $\F_{S,f}$ is
saturated however, one can get a candidate
$\Link_{S,\theta}^c(\pComp{|\Link|})$ for an associated centric
linking system by retaining information about the homotopies giving
the equivalence \mbox{$f\vert_{BP} \simeq f\vert_{BQ}\circ
B\varphi$} in the definition above. See \cite{BLO2} for details.

\begin{theorem}\cite[Theorems 7.4,7.5]{BLO2}\label{thm:p-LFGfromClSp}
For a $p$-local finite group $\plfg$, the fusion system
$\F_{S,\theta}(\pComp{|\Link|})$ is saturated and
$\Link_{S,\theta}(\pComp{|\Link|})$ is a centric linking system
associated to $\F_{S,\theta}(\pComp{|\Link|})$. Furthermore, the
\mbox{$p$-local} finite groups $(S,\F,\Link)$ and
$(S,\F_{S,\theta}(\pComp{|\Link|}),\Link_{S,\theta}^c(\pComp{|\Link|}))$
are isomorphic.
\end{theorem}

\subsection{$p$-local finite groups over abelian
groups}\label{sub:Abelian} We conclude this review by classifying
the $p$-local finite groups over an abelian $p$-group $S$. The
resulting classification shows that the strict equivalence classes
of $p$-local finite groups over $S$ are in a bijective
correspondence with the subgroups $W \leq \Aut{S}{}$ of order prime
to $p$, under the assignment
$$W \mapsto (S,\F_{S}(\WS),\Link_{S}(\WS)),$$
where $\WS$ is the semi-direct product. Here we say that two
$p$-local finite groups are strictly equivalent if they have the
same fusion system and their linking systems are isomorphic. In
particular, there are no exotic $p$-local finite groups over abelian
$p$-groups.

We begin with the following lemma which describes precisely how the
conditions in definition \ref{def:SatFS} are simplified under the
assumption that $S$ is abelian. As this result is obvious to the
experienced reader and proving it straight from the definition is an
excellent exercise for those new to $p$-local finite groups, the
proof is left to the reader.
\begin{lemma} \label{lem:Simplification}
Let $\F$ be a fusion system over a finite abelian $p$-group $S$.
Then $\F$ is saturated if and only if the following two conditions
are satisfied:
\begin{itemize}
  \item[(I${}_{ab}$)] $\Aut{S}{\F}$ has order prime to $p$.
  \item[(II${}_{ab})$] Every $\varphi \in \Hom{P}{Q}{\F}$ is the
              restriction of some $\tilde{\varphi} \in \Aut{S}{\F}$.
\end{itemize}
\end{lemma}

The following proposition follows easily.
\begin{proposition}
  If $S$ is an abelian finite $p$-group, then the assignment \mbox{$W
  \mapsto \F_S(\WS)$} gives a bijective correspondence between subgroups \mbox{$W \leq
  \Aut{S}{}$} of order prime to $p$ and saturated fusion systems over $S$.
\end{proposition}

Being the fusion system of a group, the fusion system $\F_S(\WS)$
has an obvious associated centric linking system $\Link^c_S(W \rsemi
S)$. This is in fact the only associated centric linking system and
we have the following classification result.
\begin{proposition} \label{prop:AbelianClassification}
  If $S$ is an abelian finite $p$-group, then the assignment
  $$W \mapsto \left(S,\F(W \rsemi S),\Link^c_S(W \rsemi S)\right)$$
  gives a bijective
  correspondence between subgroups \mbox{$W \leq Aut(S)$} of order prime to
  $p$ and strict equivalence classes of $p$-local finite groups over $S$. In particular, there are
  no exotic $p$-local finite groups over $S$.
\end{proposition}
\begin{proof} When $S$ is abelian, there are no proper centric
subgroups. Therefore the obstruction to uniqueness of centric
linking systems (\cite[Section 3]{BLO2}) simplifies to the
cohomology group $H^2(W;S)$. Now use a transfer argument to show
that $\Coh{W;S}$ vanishes for \mbox{$* > 0$}.
\end{proof}

\section{Retractive transfer triples} \label{sec:Setting}
In this section we introduce retractive transfer triples. First we
make precise the setting we are working in. Cohomology will always
be taken to be with $\Fp$-coefficients unless otherwise specified.
The following definition is a homotopy generalization of group
monomorphisms.
\begin{definition}
A map \mbox{$f \negmedspace : Y \to X$} between two topological
spaces is a \emph{homotopy monomorphism at $p$} if its induced map
in cohomology makes $\Coh{Y}$ a finitely generated $\Coh{X}$-module.
In the special case where $Y = BP$ is the classifying space of a
finite $p$-group, we say that $f$ is a \emph{$p$-subgroup
inclusion}.
\end{definition}
The analogy with group monomorphisms is that a group homomorphism
\mbox{$\varphi \negmedspace : P \to G$} from a finite $p$-group $P$
to a finite group $G$ is a monomorphism if and only $B\varphi$ is a
homotopy monomorphism at $p$. There are other definitions of
homotopy monomorphisms in the literature, but these are equivalent
in the setting in which we are working. Since the finite generation
hypothesis is the only one we need, we avoid complication by
considering only this definition. As the prime $p$ is fixed
throughout we will refer to these concepts simply as ``homotopy
monomorphism'' and ``subgroup inclusion.''

We will demand some additional structure on our subgroup inclusions,
namely that they allow a transfer with properties similar to that of
the transfer of a Sylow subgroup inclusion.

\begin{definition}\label{def:FrobTransfer}
Let $f \negmedspace : Y \to X$ be a map of spaces. A
\emph{retractive transfer} of $f$ is a stable map $t: \PtStable{X}
\to \PtStable{Y}$ such that $\PtStable{f} \circ t \simeq
id_{\PtStable{X}}$, and the following diagram commutes up to
homotopy
\begin{equation} \label{eq:FrobeniusDiagram}
\begin{CD}
{\PtStable{X}} @ > {\Delta} >> {\PtStable{X}\wedge\PtStable{X}} \\
@ VV {t} V @ VV {1 \wedge t} V \\
{\PtStable{Y}} @> {(f\wedge 1) \circ \Delta} >> {\PtStable{X}\wedge\PtStable{Y}.}\\
\end{CD}
\end{equation}
\end{definition}

The objects that will be the focus of our attention are defined as
follows.
\begin{definition}\label{def:RTT}
A \emph{retractive transfer triple} over a finite $p$-group $S$ is a
triple $(f,t,X)$, where $X$ is a connected, $p$-complete, nilpotent
space with finite fundamental group, $f$ is a subgroup inclusion $BS
\to X$ and $t$ is a retractive transfer of $f$.
\end{definition}

Since the space $X$ in the above definition is $p$-complete with
finite fundamental group, it follows that $\pi_1(X)$ is a finite
$p$-group \cite{DW}.

For a retractive transfer triple $(f,t,X)$ over a finite $p$-group
$S$, we ask the following questions:
\begin{itemize}
  \item[$\bullet$] Is the fusion system $\F_{S,f}(X)$ saturated?
  \item[$\bullet$] If so, does there exist an associated centric linking system $\Link$? Is it unique?
  \item[$\bullet$] If an associated centric linking system exists, then what is the relation
  between the classifying space $\ClSp$ and $X$? Are they equivalent as
  objects under $BS$?
\end{itemize}

In the course of the following sections we will answer these
questions affirmatively in the case when $S$ is an elementary
abelian $p$-group, proving the following theorem.

\begin{theorem} \label{thm:ElAb}
Let $(f,t,X)$ be a retractive transfer triple over be a finite
elementary abelian $p$-group $V$, and put
 \[W := \Aut{V}{\F_{V,f}(X)} = \{ \varphi \in \Aut{V}{} \mid f\circ \varphi \simeq f \} \]
Then the following hold
\begin{itemize}
  \item[(i)]
$W$ has order prime to $p$.
  \item[(ii)]
$\F_{V,f}(X)$ is equal to the saturated fusion system \mbox{$\F_V(W
\rsemi V)$}
  \item[(iii)]
$\F_{V,f}(X)$ has an associated centric linking system, which is
unique up to isomorphism, with classifying space \mbox{$\pComp{B(W
\rsemi V)}$}.
  \item[(iv)]
There is a homotopy equivalence \mbox{$h \colon \pComp{B(W \rsemi
V)} \xrightarrow{\simeq} X$} making the following diagram commute up
to pointed homotopy.
\[
 \xymatrix{
 & \pComp{B(W \rsemi
V)} \ar[2,0]^h_{\simeq}\\
 BV \ar[-1,1]^{Bi} \ar[1,1]^f &\\
 & X.
 }
\]
\end{itemize}
Thus the triple $(f,t,X)$ induces a \mbox{$p$-local} finite group
$(V,\F_{V,f}(X),\Link_{V,f}^c(\pComp{X}))$ over $V$ with classifying
space $X$.
\end{theorem}
\begin{proof}
The proof is by forward referencing. Lannes's theorem \cite{Lan}
shows that $W$ is the group of automorphisms of $V$ that act
trivially on $\Coh{X}$ when regarded as a subring of $\Coh{BV}$
under $f^*$. By Proposition \ref{prop:Fixed}, $W$ has order prime to
$p$, proving Part (i). By Corollary \ref{cor:CohMapExists}, $X$ has
the cohomology type of $B(W \rsemi V)$ as objects under $BV$. By
Proposition \ref{prop:MapExists} there is a map ${B(W \rsemi V)} \to
X$ realizing that cohomology isomorphism and making the uncompleted
version of the diagram in Part (iv) commute up to pointed homotopy.
This map becomes a homotopy equivalence upon $p$-completion, proving
Part (iv). Part (ii) follows directly from Parts (i) and (iv), and
Part (iii) then follows from \ref{prop:AbelianClassification}.
\end{proof}

\section{Cohomology type of retractive transfer triples}
\label{sec:Cohomology} In this section we discuss the cohomological
structure of a retractive transfer triple $(f,t,X)$ over a finite
$p$-group $S$. We first discuss general properties in
\ref{sub:CohGen}. We then specialize to the case where $S$ is
elementary abelian in \ref{sub:CohElAb} and show that in this case
$\Coh{X}$ is a ring of invariants of $\Coh{BS}$ under the action of
a group of order prime to $p$.
\smallskip

\subsection{The general case} \label{sub:CohGen}
Applying the cohomology functor $\Coh{-;\Fp}$ to
(\ref{eq:FrobeniusDiagram}) we get maps
$$\Coh{X} \xrightarrow{f^{*}} \Coh{BS} \xrightarrow{t^{*}} \Coh{X}$$
with the following properties:
\begin{itemize}
 \item[CohI\phantom{II}] $t^{*} \circ f^{*} = id$.
 \item[CohII\phantom{I}] $t^{*}$ is $\Coh{X}$-linear (Frobenius reciprocity).
 \item[CohIII] $t^{*}$ is a morphism of unstable modules over the Steenrod algebra.
 \item[CohIV] $f^*$ is a morphism of unstable algebras over the
 Steenrod algebra.
\end{itemize}
Hence $\Coh{X}$ is a direct summand of $\Coh{BS}$ as a
$\Coh{X}$-module and as a module over the Steenrod algebra. CohI
allows us to regard $\Coh{X}$ as a subring of $\Coh{BS}$ and we
will often do so without further comment.

These properties are quite restrictive and the question of which
unstable subalgebras \mbox{$R^* \subset \Coh{BS}$} over the Steenrod
algebra admit a splitting \mbox{$\Coh{BS} \to R^{*}$} as
$R^{*}$-modules and unstable modules over the Steenrod algebra is
interesting in itself. However, we focus our attention on $p$-local
finite groups.

The following finiteness properties of retractive transfer triples
will be needed later.
%
\begin{lemma}\label{lem:Fp-finite}
Let $S$ be a finite $p$-group and $(f,t,X)$ be a Frobenius
transfer triple over $S$. Then $\Coh{X}$ is Noetherian and in
particular $X$ is of finite $\Fp$-type.
\end{lemma}
\begin{proof}
By \cite[Lemma 2.6]{DW}, this follows from CohI, CohII and the
classical result that $\Coh{BS}$ is Noetherian \cite{Evens,Venkov}.
\end{proof}

\begin{lemma}\label{lem:Zp-finite}
Let $S$ be a finite $p$-group and $(f,t,X)$ be a Frobenius
transfer triple over $S$. Then $X$ is of $\mathbf{Z}_{(p)}$-finite
type.
\end{lemma}
\begin{proof}
By the universal coefficient theorem, it suffices to show that $X$
is of finite $\Fp$-type and of finite $\Q$-type. The former is Lemma
\ref{lem:Fp-finite} above. The latter is deduced in a similar way:
By a transfer argument, $BS$ has trivial $\Q$-cohomology. As in the
$\Fp$-coefficient case, $H^*(X;\Q)$ is a direct summand of
$H^*(BS;\Q)$. Hence $X$ also has trivial $\Q$-cohomology and we are
done.
\end{proof}

\subsection{The elementary abelian case} \label{sub:CohElAb} In this
subsection, we restrict ourselves to the case where $S$ is an
elementary abelian finite $p$-group $V$. In this case, we use a
theorem of Goerss-Smith-Zarati \cite{GSZ}, based on the celebrated
work of Adams-Wilkerson \cite{AW}, to prove that if $(f,t,X)$ is a
retractive transfer triple over $V$, then $\Coh{X}$ is a ring of
invariants $\Coh{BV}^{W}$ for a subgroup $W \leq \Aut{V}{}$ of order
prime to $p$. Furthermore we show that the group $W$ may be taken to
be the group of automorphisms of $V$ that act trivially on
$\Coh{X}$. We consider only the case of an odd prime. The results
still hold true at the prime $2$ and the proofs proceed in more or
less the same way, but are simpler at times. As pointed out to the
author by Nick Kuhn, these results can also be obtained, possibly
more directly, as a consequence of \cite{HLS}.

In \cite{AW}, Adams and Wilkerson study the following category.
\begin{definition}
Let $\AW$ be the category of evenly graded unstable algebras $R$
over the Steenrod algebra, that are integral domains.
\end{definition}
They also make precise the notions of ``algebraic extension'' and
``algebraic closure'' in this setting and prove the following.
\begin{proposition} \cite[Proposition 1.5]{AW}
Every object $R$ in $\AW$ has an algebraic closure $H$ in $\AW$. If
$R$ has finite transcendence degree, then so does $H$.
\end{proposition}
\begin{theorem} \cite[Theorem 1.6]{AW}
The objects $H$ in $\AW$, that are algebraically closed and of
finite transcendence degree are precisely the polynomial algebras
$\Fp[x_1,\dots,x_n]$ on generators $x_i$ of degree $2$.
\end{theorem}

Furthermore, in \cite[Theorem 1.2]{AW} they show that an algebra in
$R$ of finite transcendence degree is a ring of invariants in its
algebraic closure if and only if it satisfies two conditions, which
can be interpreted as an integral closure condition and an
inseparable closure condition. Based on that work,
Goerss-Smith-Zarati identified sufficient conditions for an unstable
algebra over the Steenrod algebra to be a ring of invariants. Before
stating their result we need some preparation.

Let $\A$ denote the mod $p$ Steenrod algebra and let $\A'$ be the
subalgebra generated by the power operations \mbox{$P^i, i \geq 0$}.
Then we have a splitting of $A'$-modules
$$\A = \A' \oplus \A'',$$
where $\A''$ is the $\Fp$-vector subspace of $\A$ generated by
admissible sequences involving the Bockstein operation. Let $\U$
denote the category of unstable $\A$-modules and $\K$ denote the
category of unstable $\A$-algebras. In both cases morphisms are of
degree zero. Let $\U'$ and $\K'$ denote the corresponding full
subcategories whose objects are evenly graded.

By \cite{LZ}, the forgetful functor \mbox{$\theta \negmedspace : \U'
\to \U$} has a right adjoint \mbox{$\tilde{\theta} \negmedspace : \U
\to \U'$}, which sends an unstable $\A$-module $M$ to the submodule
of elements $x$ of even degree satisfying \mbox{$\alpha(x) = 0$} for
all \mbox{$\alpha \in \A''$}, and morphisms to restrictions to these
submodules. For unstable $\A$-algebras the same construction gives a
right adjoint \mbox{$\tilde{\theta} \negmedspace : \K \to \K'$} to
the forgetful functor \mbox{$\theta \negmedspace : \K' \to \K$}, and
we have a commutative diagram of functors
\[
\begin{CD}
{\K'} @ > {\theta} >> {\K} @> {\tilde{\theta}} >> {\K'}\\
@ VVV @ VVV @ VVV \\
{\U'} @ > {\theta} >> {\U} @> {\tilde{\theta}} >> {\U',}\\
\end{CD}
\]
where the vertical functors are forgetful functors. As a consequence
we obtain the following lemma.
\begin{lemma} Properties \textup{CohI} to \textup{CohIV} are preserved by
$\tilde{\theta}$.
\end{lemma}
\begin{proof} This is mostly self-evident. The functor diagram is only
needed to make sense of CohII.
\end{proof}

We need to recall some things about reduced $\U$-injectives.
\begin{definition} \label{def:U-inj}
An $\A$-module $M$ is a \emph{reduced $\U$-injective} if it is an
injective object in the category $\U$, and
$$\Hom{\Sigma N}{M}{\U} = 0,$$
for every $\A$-module $N$, where $\Sigma$ denotes the suspension
functor. An unstable $\A$-algebra $R$ is a reduced $\U$-injective if
it is a reduced $\U$-injective when regarded as an $\A$-module.
\end{definition}
For an elementary abelian $p$-group $V$, the cohomology ring
$\Coh{BV}$ is a reduced \mbox{$\U$-injective} by \cite{LZ}. If
$(f,t,X)$ is a retractive transfer triple over $V$, then $\Coh{X}$
is a direct summand of $\Coh{BV}$ as $\A$-modules, and hence
$\Coh{X}$ is also a reduced $\U$-injective. This allows us to apply
the following theorem to show that $\Coh{X}$ is a ring of invariants
in $\Coh{BV}$.

\begin{theorem}\cite[Theorem 1.3]{GSZ}\label{thm:GSZ} Let $R$ be an unstable $\A$-algebra that is a
reduced $\U$-injective satisfying
\begin{itemize}
  \item[(i)] $\tilde{\theta}R$ is a Noetherian integral domain.
  \item[(ii)] $\tilde{\theta}R$ is integrally closed in its field
  of fractions.
\end{itemize}
Then there exists an integer $n$ and a subgroup \mbox{$W \leq
GL(n,\Z/p)$} such that $R$ is isomorphic to the ring of invariants
$\Coh{B(\Z/p)^n;\Fp}^W$. Furthermore, $W$ has order prime to $p$.
\end{theorem}
\begin{remark}\label{rem:n=TrDeg}
Looking closely at the proof of the theorem in \cite{GSZ} and the
tools from \cite{AW} used therein, one sees that
$\tilde{\theta}\Coh{B(\Z/p)^n;\Fp}$ is in fact the algebraic closure
of $\tilde{\theta}R^*$ in $\AW$ and $W$ is the group of
automorphisms of $(\Z/p)^n$ acting trivially on $R$. The point is
that, in their proof, Goerss-Smith-Zarati apply Theorem 1.2 of
\cite{AW}, which, as mentioned above, really gives necessary and
sufficient conditions for when the embedding of an algebra of finite
transcendence degree into its algebraic closure in $\AW$ is a Galois
extension. This is made clear in the introduction of \cite{AW}
although the authors chose to make the statement of the theorem less
technical.
\end{remark}

The following technical result is needed.
\begin{lemma}\label{lem:StillFG} Let $R$ and $H$ be unstable $\A$-algebras that are reduced
$\U$-injectives and suppose \mbox{$f \negmedspace : R \to H$} is a
morphism of $\A$-algebras making $H$ finitely generated over $R$.
Then $\tilde{\theta}f$ makes $\tilde{\theta}H$ finitely generated
over $\tilde{\theta}R$.
\end{lemma}
\begin{proof} Recall from \cite{GSZ} that there are unique
$\A'$-algebra homomorphisms
$$\pi_R \negmedspace : R \to \tilde{\theta}R~ \text{and}~ \pi_H \negmedspace : H \to \tilde{\theta}H$$
such that
$$\pi_R \circ i_R = id_{\tilde{\theta}R}~ \text{and}~ \pi_H \circ i_H = id_{\tilde{\theta}H},$$
where
$$i_R \negmedspace : \tilde{\theta}{R}~ \to R~ \text{and}~ i_H \negmedspace : \tilde{\theta}H~ \to H$$
denote the natural inclusions. (Strictly speaking this is an abuse
of notation and we should replace $\tilde{\theta}{H}$ and
$\tilde{\theta}{R}$ by $\theta\tilde{\theta}{H}$ and
$\theta\tilde{\theta}{R},$ respectively.)

Now, suppose $\{h_1,\ldots,h_n \}$ is a set of generators for $H$
over $R$. Let \mbox{$h \in \tilde{\theta}H$}. Then we can write
$$i_H(h) = \sum_{j=1}^{n} f(r_j)h_j,$$
for some \mbox{$r_j \in R$} Consequently,
\begin{align*}
  h = \pi_H \circ i_H(h) &= \pi_H \left( \sum_{j=1}^{n} f(r_j)h_j \right)\\
    &= \sum_{j=1}^{n} \pi_H(f(r_j)) \pi_H(h_j).
\end{align*}
Therefore, if we can show that \mbox{$\pi_H \circ f = \tilde{\theta
}f \circ \pi_R$}, then we can deduce that
$\{\pi_H(h_1),\ldots,\pi_H(h_n) \}$ is a set of generators for
$\tilde{\theta}H$ over $\tilde{\theta}R$ and we are done.

To prove this we first observe that by \cite[Corollary 3.3(i)]{GSZ}
there is a unique morphism of unstable $\A$-algebras \mbox{$g
\negmedspace : R \to H$} such that
$$\pi_H \circ g = \tilde{\theta}f \circ \pi_R.$$
Next we note that by construction of $\tilde{\theta}$ we have
$$i_H \circ \tilde{\theta}g = g \circ i_R,$$
from which it follows that
$$\tilde{\theta}g = \pi_H \circ i_H \circ \tilde{\theta}g = \pi_H \circ g \circ i_R = \tilde{\theta}f \circ \pi_R \circ i_R = \tilde{\theta}f.$$
But by \cite[Corollary 3.3(ii)]{GSZ} the map
$$\tilde{\theta}  : \Hom{R}{H}{\K} \longrightarrow \Hom{\tilde{\theta}R}{\tilde{\theta}H}{\K'}$$
is a bijection, so \mbox{$g = f$}, and hence
$$\pi_H \circ f = \pi_H \circ g = \tilde{\theta}f \circ \pi_R.$$
\end{proof}

\begin{proposition}\label{prop:Fixed}
Let $(f,t,X)$ be a retractive transfer triple over a finite
elementary abelian $p$-group $V$, and let \mbox{$W \leq \Aut{V}{}$}
be the subgroup of automorphisms of $V$ acting trivially on
$f^*(\Coh{X})$. The map induced by $f$ in cohomology is a split
monomorphism
 \[ f^* \colon \Coh{X} \hookrightarrow \Coh{BV}  \]
with image the ring of invariants $\Coh{BV}^W$. Furthermore, $W$ has
order prime to $p$.
\end{proposition}
\begin{proof} We already know that $f^*$ is a split monomorphism.
Let us show that $\Coh{X}$ satisfies the conditions of Theorem
\ref{thm:GSZ} above. By the remark after Definition \ref{def:U-inj},
$\Coh{X}$ is a reduced $\U$-injective. By \cite{LZ,Zar} we have
\begin{equation} \label{eq:adjointCohBS}
\tilde{\theta}\Coh{BV} \cong \Fp[x_1,\dots,x_n],
\end{equation}
where $n$ is the rank of $V$. In particular,
$\tilde{\theta}\Coh{BV}$ is a Noetherian integral domain. A similar
argument to the proof of Lemma \ref{lem:Fp-finite} shows that
$\tilde{\theta}\Coh{X}$ is also a Noetherian integral domain. It
remains only to show that $\tilde{\theta}\Coh{X}$ is integrally
closed in its field of fractions.

For this we first recall from \cite{AW}, that
$\tilde{\theta}\Coh{BV}$ is integrally closed in its field of
fractions. Now, let $x$ be in the field of fractions of
$\tilde{\theta}\Coh{X}$ and suppose that $x$ is integral over
$\tilde{\theta}\Coh{X}$. Write $x = a/b$, with $a,b \in
\tilde{\theta}\Coh{X}$. Then $\tilde{\theta}f^*(x) =
\tilde{\theta}f^*(a)/\tilde{\theta}f^*(b)$ is also integral over
$\tilde{\theta}\Coh{BV}$, and since $\tilde{\theta}\Coh{BV}$ is
integrally closed, this implies that $\tilde{\theta}f^*(x) \in
\tilde{\theta}\Coh{BV}$. We now have the equation
$\tilde{\theta}f^*(a) = \tilde{\theta}f^*(b)\tilde{\theta}f^*(x)$ in
$\tilde{\theta}\Coh{BV}$. Applying $\tilde{\theta}t^*$ and using
\mbox{$\tilde{\theta}\Coh{X}$-linearity} (CohII), we get:
 \[a = \tilde{\theta}t^*(\tilde{\theta}f^*(a)) = \tilde{\theta}t^*(\tilde{\theta}f^*(b)\tilde{\theta}f^*(x)) = b\cdot\tilde{\theta}t^*(\tilde{\theta}f^*(x)).\]
Since $\tilde{\theta}\Coh{BV}$ is an integral domain, this implies
that
 \[x = a/b = \tilde{\theta}t^*(\tilde{\theta}f^*(x)) \in \tilde{\theta}\Coh{X}.\]

Before applying Theorem \ref{thm:GSZ}, we note that since $\Coh{BV}$
is finitely generated over $\Coh{X}$, Lemma \ref{lem:StillFG}
applies to show that $\tilde{\theta}\Coh{BV}$ is finitely generated
over $\tilde{\theta}\Coh{X}$, and hence $\tilde{\theta}\Coh{BV}$ is
an algebraic extension of $\tilde{\theta}\Coh{X}$. Since
$\tilde{\theta}\Coh{BV}$ is algebraically closed in $\AW$, we
conclude that $\tilde{\theta}\Coh{BV}$ is the algebraic closure of
$\tilde{\theta}\Coh{X}$ in $\AW$.

Applying Theorem \ref{thm:GSZ} along with Remark \ref{rem:n=TrDeg}
and the observation in the preceding paragraph now completes the
proof.
\end{proof}

As an immediate corollary, a retractive transfer triple over an
elementary abelian $p$-group $V$ has the cohomology type of a
$p$-local finite group over $V$.
\begin{corollary} \label{cor:CohMapExists}
Let $(f,t,X), V$ and $W$ be as in the previous proposition. Let $G$
be the semi-direct product $G:= W \rsemi V$, and let $i$ be the
inclusion \mbox{$i \colon V \hookrightarrow G$}. There is an
isomorphism of unstable $\A$-algebras
\[ h^* \colon \Coh{X} \xrightarrow{\cong} \Coh{BG}\]
making the following diagram commute
\[
 \xymatrix{
  {\Coh{X}}  \ar@{>->}[1,1]^{f^*} \ar[2,0]^{h^*}_\cong &\\
  &          {\Coh{BV}}\\
  {\Coh{BG}} \ar@{>->}[-1,1]^{Bi^*}&&.
 }
\]
\end{corollary}
\begin{proof}
By a well known transfer argument, $Bi^*$ is a split monomorphism
with image $\Coh{BV}^W$. A $\Coh{BV}^W$-linear splitting map is
given by $\displaystyle \frac{1}{|W|}tr_V^*$, where $tr_V$ is the
transfer associated to the $|W|$-fold covering map $Bi$. For the map
$h^*$ one can take the composite $\displaystyle \frac{1}{|W|}tr_V^*
\circ f^*$.
\end{proof}

\section{Homotopy type of retractive transfer triples} \label{sec:Homotopy}
Having identified the cohomology type of a retractive transfer
triple over an elementary abelian $p$-group $V$ as that of a
$p$-local finite group over $V$ in the preceding section, in this
section we carry that result over to homotopy. More precisely, we
construct a map of spaces realizing the map $h^*$ of Corollary
\ref{cor:CohMapExists}. We follow the approach taken by
Dwyer--Miller--Wilkerson in \cite{DMW}, using Lannes technology
\cite{Lan} to pass from cohomology to homotopy.

Throughout this section, let $(f,t,X)$ be a fixed retractive
transfer triple over an elementary abelian group $V$, and let $W$ be
the group of automorphisms of $V$ acting trivially on $f^*H^*(X)
\subseteq H^*(BV)$. Recall from Proposition \ref{prop:Fixed} that
$W$ has order prime to $p$. Let $G := W \rsemi V$ be the semi-direct
product, and let $i \colon V \hookrightarrow G$ be the inclusion.
\begin{proposition}\label{prop:MapExists} There exists a map
 \mbox{$h \colon BG \to X$}
making the following diagram commute up to pointed homotopy
\[
 \xymatrix{
 & BG \ar[2,0]^h\\
 BV \ar[-1,1]^{Bi} \ar[1,1]^f &\\
 & X.
 }
\]
\end{proposition}
The rest of this section is dedicated to constructing the map $h$.

Replace \mbox{$BV \xrightarrow{f} X$} with a homotopy equivalent
fibration \mbox{$B\V \xrightarrow{\f} X$}. This can be done in such
a way (by just using the standard construction) that there is a
homotopy equivalence $s \colon BV \to B\V$ such that $f = \f \circ
s$ and a (abusively denoted) homotopy inverse $s^{-1} \colon B\V \to
BV$ such that $s^{-1} \circ s = id_{BV}$ and $s \circ s^{-1} \simeq
id_{B\V}$.

We have a fibration
\[ \Map(B\V,B\V)_{\f} \xrightarrow{\f\circ -} \Map(B\V,X)_{\f},\]
where $\Map(B\V,X)_{\f}$ is the connected component of $\Map(B\V,X)$
containing $\f$ and $\Map(B\V,B\V)_{\f}$ is the subspace of
$\Map(B\V,B\V)$ consisting of those components that map to $
\Map(B\V,X)_{\f}$. We denote the fibre over $\f$ by $\W$. This is
the space of selfmaps $g$ of $B\V$ such that $\f \circ g = \f$. Such
a map $g$ is necessarily a homotopy equivalence since $\f$ is a
homotopy monomorphism.

The next lemma can be interpreted as saying that
$BV\xrightarrow{f}X$ is centric.
\begin{lemma} \label{lem:evEqv}
For $g \in \W$, the map
\begin{equation} \label{eq:evMap}
 \Map(B\V,B\V)_g \xrightarrow{\f\circ -} \Map(B\V,X)_{\f}
\end{equation}
is a homotopy equivalence.
\end{lemma}
\begin{proof}
The map in question is adjoint to the bottom row of the commutative
diagram
\[
 \xymatrix{
 {\Map(BV,BV)_{id} \times BV} \ar[0,1]^{{\hskip .65in} ev} \ar[1,0]^{(g_*\circ c_s)\times s}_{\simeq} &   BV \ar[1,0]^{g \circ s}_{\simeq} \ar[0,1]^{Bi} \ar[1,1]^f &  BG \\
 {\Map(B\V,B\V)_{g}  \times B\V} \ar[0,1]^{{\hskip .65in} ev}  & B\V \ar[0,1]^{\f} &  X,
 }
\]
where $c_s$ is the map sending a selfmap $u \in \Map(BV,BV)_{id}$ to
its ``conjugate'' \mbox{$s \circ u \circ s^{-1} \in
\Map(B\V,B\V)_{id}$}, and $g_*$ is composition with $g$. Applying
the cohomology functor we obtain the commutative diagram
\[
  \xymatrix{
  {\Coh{X}}  \ar[0,1]^{\f^*} \ar[1,0]^{h^*}_{\cong} \ar[1,1]^{f^*}   & {\Coh{B\V}} \ar[1,0]^{(g\circ s)^*}_{\cong} \ar[0,1]^{ev^*{\hskip 0.8 in}}      & {\Coh{\Map(B\V,B\V)_{g}} \otimes \Coh{B\V}} \ar[1,0]^{(g_*\circ c_s)^*\otimes s^*}_{\cong}\\
  {\Coh{BG}} \ar[0,1]^{Bi^*}    & {\Coh{BV}}  \ar[0,1]^{ev^*{\hskip 0.8 in}}    & {\Coh{\Map(BV,BV)_{id}} \otimes \Coh{BV}}, \\
 }
\]
where $h^*$ is the isomorphism from Corollary
\ref{cor:CohMapExists}. Taking adjoints and restricting to
components we obtain the commutative diagram

\begin{equation} \label{eq:TVDiagram} \raisebox{0.90in}{
  \xymatrix{
  {T^{\V}_{\f}(\Coh{X})} \ar[0,2]^{T^{\V}_{g}(\f^*)} \ar[1,0]^{\eta}_{\cong}  &  & {T^{\V}_{g}(\Coh{B\V})} \ar[1,0]^{\eta}_{\cong} \ar[0,2]^{\lambda_{g} {\hskip 0.2 in}}    & & {\Coh{\Map(B\V,B\V)_{g}}} \ar[3,0]^{(g_*\circ c_s)^*}_{\cong}
  \\
  {T^V_{f}(\Coh{X})} \ar[0,2]^{T^V_{g\circ s}(\f^*)} \ar[2,0]^{T^V_{Bi}(h^*)}_{\cong} \ar[2,2]^{T^V_{id}(f^*)}  &  & {T^V_{g\circ s}(\Coh{B\V})} \ar[2,0]^{T^V_{id}((g\circ s)^*)}_{\cong}
  \\
  \\
  {T^V_{Bi}(\Coh{BG})} \ar[0,2]^{T^V_{id}(Bi^*)}    & & {T^V_{id}(\Coh{BV})}  \ar[0,2]^{\lambda_{id} {\hskip 0.2 in}}    &  & {\Coh{\Map(BV,BV)_{id}}},
 }}
\end{equation}
where the Lannes functors $T^V$ and $T^{\V}$ are the left adjoints
to the functors
 $-\otimes \Coh{BV}$ and $-\otimes \Coh{B\V}$, respectively,
and $\eta$ is the natural isomorphism of functors induced by the
isomorphism $s^*$. The subscript notation $T^V_{\f}$ denotes the
component of $T^V$ corresponding to the map $\f^*$, and the maps
$\lambda_{id}$ and $\lambda_{g}$ are appropriate restrictions of the
adjoint to the evaluation map. The reader is referred to Section 3
of \cite{DMW} for further explanation of the notation, to
\cite{Schwartz} for a good exposition on Lannes functor technology
and to \cite{Lan} for original paper by Lannes.

By Lannes's comparison theorem, \cite[Theorem 3.3.2]{Lan}, the map
(\ref{eq:evMap}) is a homotopy equivalence if and only if the
composite of the top row in (\ref{eq:TVDiagram}) is an isomorphism.
(We have used the fact that the spaces $X$ and $\Map(B\V,B\V)_g
\simeq BV$ are already $p$-complete here.) The proof is completed by
recalling from \cite{Lan} that the maps in the bottom row are both
isomorphisms.
\end{proof}

The next proposition corresponds to Theorem 2.9 of \cite{DMW}.
\begin{proposition} \label{prop:WfisW}
The space $\W$ is homotopy discrete, and the cohomology functor
induces a natural monomorphism
 \[ \Coh{-} \colon \pi_0\W \longrightarrow \Aut{\Coh{B\V}}{\K} \]
with image $W$ under the identification
 \[  \Aut{\Coh{B\V}}{\K} \xrightarrow{\cong} \Aut{\Coh{BV}}{\K} \cong \Aut{V}{} \]
induced by $s \colon BV \xrightarrow{\simeq} B\V$.
\end{proposition}
\begin{proof}
The first claim follows from Lemma \ref{lem:evEqv}. Also from Lemma
\ref{lem:evEqv} and the long exact sequence in homotopy induced by
the fibre sequence
\[\W \xrightarrow{~ incl~}  \Map(B\V,B\V)_{\f} \xrightarrow{~\f\circ -~} \Map(B\V,X)_{\f} \]
we deduce that the map
\[ \pi_0\W \xrightarrow{~ incl~} \pi_0\Map(B\V,B\V)_{\f} \]
is a bijection. In the commutative diagram
\[
 \xymatrix{
 [B\V,B\V] \ar[0,2]^{\f\circ- {\hskip .35in}} \ar[1,0]^{\Coh{-}} && [B\V,X]\ar[1,0]^{\Coh{-}}\\
 \Aut{\Coh{B\V}}{\K} \ar[0,2]^{-\circ \f^* {\hskip 0.35in}} && \Hom{\Coh{X}}{\Coh{B\V}}{\K}
 }
\]
the left vertical arrow is an isomorphism by Miller's theorem
\cite{Mil} and the right vertical arrow is an isomorphism by
Lannes's theorem \cite{Lan}. Consequently a map $g\colon B\V \to
B\V$ belongs to $\Map(B\V,B\V)_{\f}$ if and only if $g^*\circ \f^* =
\f^*$ in cohomology, or equivalently, if and only if $g^* \in W$
(under the identification above).
\end{proof}

Observing that $\W$ is a grouplike topological monoid, we have a
contractible CW-complex $E\W$ on which $\W$ acts freely, which
allows us to form the classifying space \mbox{$B\W = E\W/\W$} and
the Borel construction $E\W \times_{\W} B\V$. By construction $\f$
induces a map $\tilde{h} \colon E\W \times_{\W} B\V \to X $ fitting
into a commutative diagram
\[
 \xymatrix{
 & E\W \times_{\W} B\V \ar[2,0]^{\tilde{h}}\\
 B\V \ar[-1,1]^{B\iota} \ar[1,1]^f &\\
 & X,
 }
\]
where $B\iota$ is the obvious map. Proposition \ref{prop:WfisW}
implies that $E\W \times_{\W} B\V$ is homotopy equivalent to $BG$.
Using the fact that these are classifying spaces of finite groups,
one sees that this homotopy equivalence can in fact be realized by a
map \mbox{$k \colon BG \to E\W \times_{\W} B\V $} making the top
rectangle in the following diagram commute up to pointed homotopy
\[
 \xymatrix{
 {BV} \ar[0,2]^{Bi} \ar[1,0]^s_{\simeq} &&{BG} \ar[1,0]^k_{\simeq}\\
 {B\V} \ar[0,2]^{B\iota} \ar[1,2]^{\f} &&E\W \times_{\W} B\V \ar[1,0]^{\tilde{h}}\\
                              &&X.
 }
\]
Since $\f \circ s = f$, the composite $\tilde{h} \circ k$ gives the
desired map $h$ in Proposition \ref{prop:MapExists}.

\section{$p$-local finite groups induce retractive transfer triples} \label{sec:p-lfg-to-RTT}
In Section \ref{sec:Setting} the notion of a retractive transfer
triple over a finite $p$-group $S$ was introduced and in Sections
\ref{sec:Cohomology} and \ref{sec:Homotopy} it was shown that, in
the case where $S$ is elementary abelian, such a triple induces a
$p$-local finite group. In this section we consider the reverse
implication and prove the following theorem.
\begin{theorem} \label{thm:p-lfg-to-RTT}
Let $\plfg$ be a $p$-local finite group. Then the natural inclusion
\mbox{$\theta \negmedspace : BS \to \pComp{|\Link|}$} has a
retractive transfer $t$, and $(\theta,t,\pComp{|\Link|})$ is a
retractive transfer triple.
\end{theorem}

There are two parts to the proof. In Subsection \ref{sub:Technical},
which is joint work with Ran Levi, we show that $\ClSp$ and the
inclusion \mbox{$\theta \negmedspace : BS \to \pComp{|\Link|}$} of
the Sylow subgroup satisfy the technical conditions of a retractive
transfer triple. Most notably we show that the classifying space of
a $p$-local finite group is both a torsion space and a nilpotent
space, a result which is of independent interest.

In Subsection \ref{sub:TransferConstruction} we apply results from
\cite{KR:ClSpec} to obtain a stable retraction $t$ of the inclusion
\mbox{$\theta \negmedspace : BS \to \pComp{|\Link|}$}, and show that
it satisfies Frobenius reciprocity. These results combine to
complete the proof of Theorem \ref{thm:p-lfg-to-RTT}

\subsection{Technical conditions} \label{sub:Technical}
Let $(S,\F,\Link)$ be a $p$-local finite group. In this subsection,
which is joint work with Ran Levi, we verify that the space
$\pComp{|\Link|}$ and the natural map
$$\theta \negmedspace : BS \to \pComp{|\Link|}$$
satisfy the technical conditions of retractive transfer triples. It
has already been shown in \cite{BLO2} that $\pComp{|\Link|}$ is
$p$-complete, that the fundamental group of $\pComp{|\Link|}$ is
finite and that $\theta$ is a homotopy monomorphism. We proceed to
show that $\pComp{|\Link|}$ is nilpotent.

\begin{lemma}\label{lem:ClSpFiniteHomology}
For a $p$-local finite group $\plfg$, the groups $H^k(\ClSp;\Z)$ and
$H_k(\ClSp;\Z)$ are finite $p$-groups for all \mbox{$k \geq 1$}.
\end{lemma}
\begin{proof} As is shown in \cite{BLO2}, $\ClSp$ is a stable retract of $BS$. In
particular $\Coh{\ClSp;\Z}$ is a subring of $\Coh{BS;\Z}$. Since
$H^k(BS;\Z)$ is a finite  $p$-group for \mbox{$k \geq 1$} (a
transfer argument shows that $H^k(BS;\Z)$ is $p$-torsion and finite
generation is evident from the cell structure of $BS$,) it follows
that $H^k(\ClSp;\Z)$ is a finite  $p$-group for \mbox{$k \geq 1$}.

The same argument works in homology.
\end{proof}

\begin{proposition} \label{prop:ClSpTorsion}
For a $p$-local finite group $\plfg$, the homotopy
groups $\pi_k(\ClSp)$ are finite $p$-groups for all \mbox{$k \geq
1$}. In particular, $\ClSp$ is a torsion space.
\end{proposition}
\begin{proof}
We first reduce to the case where $\ClSp$ is simply connected. It is
shown in \cite{BLO2} that $\pi_1(X)$ is a finite $p$-group. Letting
$\widetilde{X}$ be a universal cover of $\ClSp$, it therefore
suffices to show that the homotopy groups of $\widetilde{X}$ are all
finite $p$-groups. But it is shown in \cite{BCGLO2} that
$\widetilde{X}$ is again the classifying space of a $p$-local finite
group, so we can reduce to the simply connected case.

Assume therefore that $\ClSp$ is simply connected. By Lemma
\ref{lem:ClSpFiniteHomology}, $H_k(\ClSp;\Z)$ is a finite $p$-group
for all \mbox{$k \geq 1$}. Since $\ClSp$ is simply connected, we can
apply the Hurewicz theorem modulo the class of finite abelian
$p$-groups, and deduce that the homotopy groups of $\ClSp$ are all
finite $p$-groups.
\end{proof}

\begin{corollary}
The classifying space of a $p$-local finite group is nilpotent.
\end{corollary}
\begin{proof}
This follows from Proposition \ref{prop:ClSpTorsion}, the fact that
any finite $p$-group is nilpotent and the fact that the action of
any finite $p$-group on a finite abelian $p$-group is nilpotent.
\end{proof}

\subsection{The retractive transfer of a $p$-local finite group}
\label{sub:TransferConstruction} In this subsection we show that the
natural inclusion of a Sylow subgroup into the classifying space of
a $p$-local finite group has a retractive transfer. We use results
from \cite{KR:ClSpec}, which were actually originally developed for
this purpose but have turned out to be perhaps more interesting than
their intended goal, and are therefore published separately.

\begin{remark}
There is a slight difference between the spectra appearing in this
paper and the ones appearing in \cite{KR:ClSpec}. The difference
arises because here we add a basepoint to our spaces before forming
suspension spectra. The effect at the level of stable homotopy is to
add a sphere wedge summand to all spectra in sight. It is easy to
check that all the results quoted from \cite{KR:ClSpec} carry over
in the form stated in this section.
\end{remark}

 Let $\plfg$ be a $p$-local finite group, and let
$\StableCharIdem$ be the idempotent of $\PtStable{BS}$ induced by
the characteristic idempotent of $\F$, as defined in
\cite{KR:ClSpec}. We refer to $\StableCharIdem$ as the \emph{pointed
stable idempotent of $\F$}. It has the following properties, which
determine $\StableCharIdem$ uniquely:
\begin{itemize}
  \item[(a)] $\StableCharIdem$ is a $\Zp$-linear combination of homotopy
classes of maps of the form \mbox{$\PtStable{B\varphi} \circ tr_P$},
where $P$ is a nontrivial subgroup of $S$, \mbox{$\varphi \in
\HomF{P}{S}$} and $tr_P$ denotes the transfer of the inclusion
\mbox{$P\leq S$}.
  \item[(b)]
For each subgroup \mbox{$P \leq S$} and each \mbox{$\varphi \in
\HomF{P}{S}$}, the restrictions
\mbox{$\StableCharIdem\circ{\PtStable{Bi_P}}$} and
\mbox{$\StableCharIdem\circ{\PtStable{B\varphi}}$} are homotopic as
maps \mbox{$\PtStable{BP} \to \PtStable{BS}$}.
 \item[(c)]
$\StableCharIdem$ has augmentation 1.
\end{itemize}
The augmentation in Property (c) corresponds to an augmentation of
the double Burnside ring. As Property (c) is not used in this paper,
the reader is referred to \cite{KR:ClSpec} for the details. We refer
to Property (b) as \emph{$\F$-stability}.

The \emph{pointed classifying spectrum of $\F$} is the stable
summand $\PtClSpec{\F}$ of $\PtStable{BS}$ induced by
$\StableCharIdem$. This is the infinite mapping telescope of
$\StableCharIdem$:
$$\PtClSpec{\F} = \text{HoColim} \left( \PtStable{BS} \xrightarrow{\StableCharIdem} \PtStable{BS} \xrightarrow{\StableCharIdem} \PtStable{BS} \xrightarrow{\StableCharIdem} \cdots \right),$$
and as such it comes equipped with a \emph{pointed structure map},
\mbox{$\PtStable{BS} \xrightarrow{\sigma} \PtClSpec{\F},$} which is
the structure map of the homotopy colimit, and a unique (up to
homotopy) map \mbox{$\PtClSpec{\F} \xrightarrow{t} \PtStable{BS} $}
such that \mbox{$\sigma \circ t \simeq id_{\PtClSpec{\F}}$} and
\mbox{$t \circ \sigma \simeq \StableCharIdem.$}

Just as in \cite{KR:ClSpec}, one can show that the pointed structure
map \mbox{$\PtStable{BS} \xrightarrow{\sigma} \PtClSpec{\F}$} is
equivalent to the infinite pointed suspension \mbox{$\PtStable{BS}
\xrightarrow{\PtStable{\theta}} \PtStable{\ClSp}$} of the natural
inclusion \mbox{$BS \xrightarrow{\theta} \ClSp$}, as objects under
$\PtStable{BS}$. We may therefore replace \mbox{$\PtStable{BS}
\xrightarrow{\sigma} \PtClSpec{\F}$} by \mbox{$\PtStable{BS}
\xrightarrow{\PtStable{\theta}} \PtStable{\ClSp}$} in the discussion
above, and obtain a unique (up to homotopy) map
$$t : \PtStable{\ClSp} \longrightarrow \PtStable{BS},$$
such that
$$\PtStable{\theta} \circ t \simeq id_{\PtStable{\ClSp}} \hspace{1cm} \text{and} \hspace{1cm} t \circ \PtStable{\theta} = \StableCharIdem.$$
We proceed to show that $t$ satisfies the Frobenius reciprocity
relation illustrated in (\ref{eq:FrobeniusDiagram}) and thus $t$ is
a retractive transfer for $\theta$.

\begin{proposition}\label{prop:omegaFrob}
The idempotent $\StableCharIdem$ satisfies the Frobenius reciprocity
relation
$$(\StableCharIdem \wedge \StableCharIdem)\circ\Delta \simeq (\StableCharIdem\wedge{1})\circ\Delta\circ\StableCharIdem,$$
where \mbox{$\Delta\negmedspace : \PtStable{BS} \to
\PtStable{BS}\wedge\PtStable{BS}$} is the image of the diagonal of
$BS$ under the infinite suspension functor $\PtStable$.
\end{proposition}
\begin{proof}
Recall that we can write $\StableCharIdem$ as a linear combination
with $\Zp$-coefficients of maps \mbox{$\PtStable{B\varphi} \circ
tr_P \in \stablemaps{BS}{BS}$}, where \mbox{$P \leq S$} and
\mbox{$\varphi \in \HomF{P}{S}$}. For such a map
\mbox{$\PtStable{B\varphi} \circ tr_P$} we have, by $\F$-stability
of $\StableCharIdem$,
\begin{equation} \label{eq:omegaInvariant}
\StableCharIdem \circ{\PtStable{B\varphi}} \simeq \StableCharIdem
\circ \PtStable{Bi_P},
\end{equation}
where $i_P$ is the inclusion $P \leq S$. We will take advantage of
this and the fact \cite{Ad} that the transfer $tr_P$ of the
inclusion $i_P$ satisfies the Frobenius relation
\begin{equation} \label{eq:FrobClassic}
(1 \wedge tr_P) \circ \Delta_S \simeq (\PtStable{Bi_P} \wedge 1)
\circ \Delta_P \circ tr_P,
\end{equation}
where $\Delta_P$ and $\Delta_S$ are the diagonals of $\PtStable{BP}$
and $\PtStable{BS}$ respectively. We will also use the fact that,
since $\PtStable {B\varphi}$ has a desuspension, it commutes with
the diagonals as follows
\begin{equation} \label{eq:phiCommutes}
\Delta_S \circ \PtStable{B\varphi} \simeq (\PtStable{B\varphi}
\wedge \PtStable{B\varphi}) \circ \Delta_P.
\end{equation}

Now,
\begin{align*}
  (\StableCharIdem\wedge1) \circ \Delta_S \circ \PtStable{B\varphi} \circ tr_P
  &\stackrel{(\ref{eq:phiCommutes})}{\simeq}(\StableCharIdem\wedge1)\circ(\PtStable{B\varphi}\wedge \PtStable{B\varphi}) \circ \Delta_P \circ tr_P\\
  &\stackrel{\phantom{(\ref{eq:phiCommutes})}}{\simeq} ((\StableCharIdem\circ \PtStable{B\varphi})\wedge \PtStable{B\varphi}) \circ \Delta_P \circ tr_P\\
  &\stackrel{(\ref{eq:omegaInvariant})}{\simeq} ((\StableCharIdem\circ \PtStable{Bi_P}) \wedge \PtStable{B\varphi}) \circ \Delta_P \circ tr_P\\
  &\stackrel{\phantom{(\ref{eq:phiCommutes})}}{\simeq} (\StableCharIdem\wedge \PtStable{B\varphi}) \circ (\PtStable{Bi_P}\wedge1) \circ \Delta_P\circ tr_P\\
  &\stackrel{(\ref{eq:FrobClassic})}{\simeq} (\StableCharIdem\wedge \PtStable{B\varphi}) \circ (1\wedge tr_P) \circ \Delta_S\\
  &\stackrel{\phantom{(\ref{eq:phiCommutes})}}{\simeq} (\StableCharIdem \wedge (\PtStable{B\varphi} \circ tr_P)) \circ \Delta_S.
\end{align*}
By summing over the different $\PtStable{B\varphi} \circ tr_P$, we
get the desired result.
\end{proof}

\begin{corollary}
Let $\plfg$ be a $p$-local finite group and let \mbox{$t
\negmedspace : \PtStable{\ClSp} \to \PtStable{BS}$} be as
constructed above. Then $t$ is a retractive transfer of the
inclusion \mbox{$\theta \negmedspace : BS \longrightarrow \ClSp.$}
\end{corollary}
\begin{proof}
We deduce the Frobenius reciprocity relation
 \[(1\wedge t) \circ \Delta_{\ClSp} \simeq (\PtStable{\theta}\wedge 1) \circ \Delta_S \circ t\]
from
\begin{equation} \label{eq:IdemFR}
(\StableCharIdem\wedge\StableCharIdem)\circ\Delta_S \simeq
(\StableCharIdem\wedge{1})\circ\Delta_S\circ\StableCharIdem
\end{equation}
as follows. Applying \mbox{$(\PtStable{\theta} \wedge 1) \circ -
\circ t$} to the left hand side of (\ref{eq:IdemFR}) and rewriting,
we get
\begin{align*} 
  (\PtStable{\theta} \wedge 1) \circ (\StableCharIdem\wedge\StableCharIdem)\circ\Delta_S \circ t
  &\simeq (\PtStable{\theta} \wedge 1) \circ ( (t \circ \PtStable{\theta}) \wedge (t \circ \PtStable{\theta}) )\circ\Delta_S \circ t\\
  &\simeq ((\PtStable{\theta} \circ t)\wedge t) \circ (\PtStable{\theta}\wedge\PtStable{\theta})\circ\Delta_S\circ t\\
  &\simeq (1\wedge t) \circ \Delta_{\ClSp} \circ \PtStable{\theta} \circ t \\
  &\simeq (1\wedge t) \circ \Delta_{\ClSp}.
\end{align*}
Doing the same with the right hand side yields
\begin{align*}\label{eq:RHS}
(\PtStable{\theta} \wedge 1) \circ
(\StableCharIdem\wedge{1})\circ\Delta_S\circ\StableCharIdem \circ t
&\simeq
((\PtStable{\theta} \circ t \circ \PtStable{\theta})\wedge 1) \circ\Delta_S\circ (t \circ \PtStable{\theta} \circ t)\\
&\simeq (\PtStable{\theta}\wedge 1) \circ\Delta_S\circ t.\\
\end{align*}
Collecting these equivalences, we have
 \[(1\wedge t) \circ \Delta_{\ClSp} \simeq (\PtStable{\theta} \wedge 1) \circ \Delta_S \circ t. \]

\end{proof}

\end{document}